\input amstex\documentstyle{amsppt}  
\pagewidth{12.5cm}\pageheight{19cm}\magnification\magstep1
\topmatter
\title A survey of total positivity\endtitle
\author G. Lusztig\endauthor
\address{Department of Mathematics, M.I.T., Cambridge, MA 02139}\endaddress
\endtopmatter   
\document
\define\we{\wedge}

\define\mpb{\medpagebreak}

\define\Lie{\text{\rm Lie }}

\define\wt{\widetilde}
\define\sqc{\sqcup}

\define\lb{\linebreak}

\define\part{\partial}

\define\ra{\rangle}

\define\m{\mapsto}
\define\do{\dots}
\define\la{\langle}
\define\bsl{\backslash}

\define\sub{\subset}    

\define\T{\times}
\define\ti{\tilde}
\define\nl{\newline}
\redefine\i{^{-1}}

\define\Ad{\text{\rm Ad}}
\define\Hom{\text{\rm Hom}}

\define\a{\alpha}

\redefine\c{\chi}
\define\g{\gamma}
\redefine\d{\delta}

\define\p{\pi}

\define\s{\sigma}
\redefine\t{\tau}

\redefine\l{\lambda}

\define\x{\xi}

\define\vt{\vartheta}

\redefine\L{\Lambda}
\define\Ph{\Phi}

\define\ii{\bold i}

\define\CC{\bold C}

\define\RR{\bold R}

\define\cb{\Cal B}
\define\cc{\Cal C}

\define\tu{\ti u}

\define\tA{\ti A}

\define\che{\check}

\define\FG{FG}
\define\GK{GK}
\define\HI{H}
\define\LU{L1}
\define\CB{L2}
\define\KA{Ka}
\define\LO{Lo}
\define\PE{P}
\define\SC{S}
\define\WH{W}

The theory of totally positive matrices originated in the 1930's in the work of
I.J.Schoenberg and Gantmacher-Krein. Since then this theory has found applications to such
diverse areas as statistics, game theory, mathematical economics, stochastic processes (see
Karlin \cite{\KA}). This survey covers the following topics:

(i) an exposition of the early work of Schoenberg and Gantmacher-Krein and the work of 
A.Whitney on totally positive matrices (see \S1);

(ii) a generalization of the known reduction to a canonical form of a symplectic 
nondegenerate bilinear form to non-symplectic bilinear forms based on ideas from total 
positivity (see \S2);

(iii) an exposition of results in \cite{\LU} which extend the notion of total positivity to
a Lie group context based on the theory of quantum groups and canonical bases (see \S3);

(iv) a generalization of (ii) placing it in a Lie group context (see \S5);

(v) an exposition of recent results of Fock and Goncharov \cite{\FG} which provide a
beautiful application of the theory of total positivity in Lie groups to the study of 
homomorphisms of the fundamental group of a closed surface into a Lie group (see \S6, \S7).

\subhead 1\endsubhead
Let $V$ be a real vector space with a totally ordered basis $e_1,e_2,\do,e_n$. For any 
$k\in[1,n]$ the $k$-th exterior power $\L^kV$ has a basis 
$(e_{r_1}\we e_{r_2}\we\do\we e_{r_k})$ indexed by the sequences $r_1<r_2<\do<r_k$ in 
$[1,n]$. Let $G=GL(V)$. Let $G_{\ge0}$ (resp. $G_{>0}$) be the set of all $A\in G$ such 
that for any $k\in[1,n]$ the coefficients of $\L^kA:\L^kV@>>>\L^kV$ with respect to the 
basis above are $\ge0$ (resp. $>0$); one then also says that $A$ is totally $\ge0$ (resp. 
totally $>0$). This definition was given by I. J. Schoenberg \cite{\SC} in 1930 in 
connection with his solution of a problem of P\'olya. For $v\in V$, $v=\sum_rv_re_r$, 
$v_r\in\RR$, let $\vt_v$ be the largest integer $k\ge0$ for which there exists a sequence 
$r_1<r_2<\do<r_k$ in $[1,n]$ such that $v_{r_1}v_{r_2}<0,\do,v_{r_{k-1}}v_{r_k}<0$; in 
other words, $\vt_v$ is the number of sign changes in the sequence $v_1,v_2,\do,v_n$ with
$0$'s removed.
P\'olya asked for a characterization of those $A\in G$ which diminish $\vt_v$ that is, are
such that $\vt_{A(v)}\le\vt_v$ for any $v\in V$. Schoenberg showed that these $A$ are 
exactly those such that for any $k\in[1,n]$, $\L^kA:\L^kV@>>>\L^kV$ does not have two 
coefficients of opposite signs with respect to the basis above (in particular any 
$A\in G_{\ge0}$ has the required diminishing property). 

A second source of the idea of totally positivity appeared in the work of Gantmacher and 
Krein \cite{\GK} in 1935. Motivated by the study of vibrations of a mechanical system, 
Gantmacher and Krein were looking for a condition on $A\in G$ which guarantees that 

(a) $A$ has distinct real positive eigenvalues;
\nl
and

-is preserved by a small perturbation of $A$;

-can be easily tested in terms of the coefficients of $A$.
\nl
They showed that the condition that $A\in G_{>0}$ has the required properties.

We now sketch a proof of the Gantmacher-Krein theorem asserting that if $A\in G_{>0}$ then
(a) holds. Let $c_1,c_2,\do,c_n$ be the eigenvalues of $A$ arranged so that 
$|c_1|\ge|c_2|\ge\do\ge|c_n|$. For any $k\in[1,n]$ the eigenvalues of $\L^kA$ are
$c_{r_1}c_{r_2}\do c_{r_k}$ for various $r_1<r_2<\do<r_k$ in $[1,n]$; if $k<n$, the first 
two eigenvalues of $\L^kA$ in the decreasing order of absolute value are 

$c_1c_2\do c_k$, $c_1c_2\do c_{k-1}c_{k+1}$.
\nl
By Perron's theorem \cite{\PE} on square matrices with all entries in $\RR_{>0}$ applied to
$\L^kA$, we have $c_1c_2\do c_k\in\RR_{>0}$ and 
$c_1c_2\do c_{k-1}c_k>|c_1c_2\do c_{k-1}c_{k+1}|$. Using this we see by induction on $k$ 
that all $c_r$ are in $\RR_{>0}$ and that $c_1>c_2>\do>c_n$.

\mpb

Let $I=[1,n-1]$. For $i\in I$, $a\in\RR$ we define $x_i(a)$ (resp. $y_i(a)$) in $G$ by
$x_i(a)e_i=e_i+ae_{i+1},x_i(a)e_h=e_h$ for $h\ne i$ (resp. by $y_i(a)e_{i+1}=e_{i+1}+ae_i$,
$y_i(a)e_h=e_h$ for $h\ne i+1$). Let $T$ (resp. $T_{>0}$) be the group of all $A\in G$ such
that $Ae_i=a_ie_i$ for all $i$ where $a_i\in\RR-\{0\}$ (resp. $a_i\in\RR_{>0}$).

In 1952, Ann Whitney \cite{\WH} (a Ph.D. student of Schoenberg) proved that $G_{>0}$ is 
dense in $G_{\ge0}$ using the statements (b),(c) below.

(b) {\it $G_{\ge0}$ is the submonoid with $1$ of $G$ generated by $x_i(a),y_i(a)$ with 
$i\in I$, $a\in\RR_{\ge0}$ and by $T_{>0}$;}

(c) {\it $G_{>0}$ consists of all products}
$$x_{i_1}(a_1)x_{i_2}(a_2)\do x_{i_N}(a_N)ty_{i_1}(b_1)y_{i_2}(b_2)\do y_{i_N}(b_N)$$
{\it where $i_1,i_2,\do,i_N$ is the sequence $1,2,\do,n-1,1,2,\do,n-2,\do,1,2,1$; 
$N=n(n-1)/2$; $a_s>0,b_s>0$ for all $s$ and $t\in T_{>0}$.}
\nl
(Actually (c) is only implicit in \cite{\WH}; it is stated explicitly in Loewner
\cite{\LO}.) Note that 

(d) {\it in (c) we can replace $1,2,\do,n-1,1,2,\do,n-2,\do,1,2,1$ by 
$n-1,n-2,\do,1,n-1,n-2,\do,2,\do,n-1,n-2,n-1$ and we get a true statement.}
\nl
Indeed if we define $G_{>0}$ in terms of the reverse order $e_n,e_{n-1},\do,e_1$ of the 
basis we get the same $G_{>0}$ as before.

\subhead 2\endsubhead
We preserve the notation of \S1. Define an involution $i\m i^*$ of $[1,n]$ by $i^*=n+1-i$.
Let $\la,\ra:V\T V@>>>\RR$ be a bilinear form. When $\la,\ra$ is symplectic, nondegenerate,
a standard result shows that there exists a basis $(v_r)_{r\in[1,n]}$ of $V$ such that 
$\la v_r,v_s\ra=0$ if $r\ne s^*$ and $\la v_r,v_{r^*}\ra=(-1)^r$ for all $r$. We want to 
prove a similar result without the assumption that $\la,\ra$ is symplectic. Instead we 
shall assume that $\la,\ra$ is totally $>0$ in the following sense: for any $k\in[1,n]$ and
any two sequences $r_1<r_2<\do<r_k$, $s_1<s_2<\do<s_k$ in $[1,n]$, the determinant of the 
matrix 
$$((-1)^{r_m}\la e_{r_m^*},e_{s_{m'}}\ra)_{m,m'\in[1,k]}$$
is $>0$. Note that the totally $>0$ bilinear forms form a nonempty open set in the space of
nondegenerate bilinear forms on $V$. 

We now fix a totally $>0$ bilinear form $\la,\ra$. We prove for it the following result.

(a) {\it There exists a basis $v_1,v_2,\do,v_n$ of $V$ such that $\la v_r,v_s\ra=0$ if 
$r\ne s^*$, $\la v_r,v_{r^*}\ra=(-1)^rz_r\in\RR-\{0\}$ for all $r$ and}
$$0<z_1z_n\i<z_2z_{n-1}\i<\do<z_nz_1\i.$$
If $X:E@>>>E'$ is an isomorphism of finite dimensional vector spaces let
$\che{X}:E^*@>>>E'{}^*$ be its transpose inverse. 

There is a unique isomorphism $C:V^*@>>>V$ such that $\la v,v'\ra=C\i(v)(v')$ for all
$v,v'$ in $V$. Let $e'_1,e'_2,\do,e'_n$ be the basis of $V^*$ dual to the basis 
$e_1,e_,\do,e_n$ of $V$. Define an isomorphism $C_0:V^*@>>>V$ by $C_0(e'_r)=(-1)^re_{r^*}$.
We have $\che{C_0}=(-1)^{n+1}C_0\i$. Define an involution $A\m\tA$ of $G$ by 
$\tA=C_0\che{A}C_0\i$. For $1\le i\le n-1,a\in\RR$ we have $\wt{x_i(a)}=x_{n-i}(a)$,
$\wt{y_i(a)}=y_{n-i}(a)$; moreover $\ti t\in T_{>0}$ for $t\in T_{>0}$. Using this and 1(d)
we see that $A\m\tA$ maps $G_{>0}$ into itself.

Define $A\in G$ by $C=C_0\che{A}$. We can write 
$A(e_r)=\sum_{r,s}a_{rs}e_s$ with $a_{rs}\in\RR$. Then
$\che{A}\i(e'_r)=\sum_{r,s}a_{sr}e'_s$. For $r,s$ in $[1,n]$ we have
$$(-1)^r\la e_{r^*},e_s\ra=(-1)^r\che{A}\i C_0\i(e_{r^*})(e_s)=\che{A}\i(e'_r)(e_s)=
\sum_ha_{hr}e'_h(e_s)=a_{sr}.$$
This and our hypothesis on $\la,\ra$ shows that $A\in G_{>0}$. We have 
$$C\che{C}=(-1)^{n+1}C_0\che{A}C_0\i A.$$
Since $A$ and $C_0\che{A}C_0\i$ are in $G_{>0}$, so is their product. Thus, 
$(-1)^{n+1}C\che{C}\in G_{>0}$. By the Gantmacher-Krein theorem for 
$(-1)^{n+1}C\che{C}\in G_{>0}$ there exists a basis $v_1,v_2,\do,v_n$ of $V$ and real 
numbers $c_1>c_2>\do>c_n>0$ such that $(-1)^{n+1}C\che{C}v_r=c_rv_r$ for all $r$.
Define $v'_r\in V^*$ by $v'_r(v_s)=\d_{rs}$. We have  
$\che{C}Cv'_r=(-1)^{n+1}c_r\i v'_r$. Let $v''_r=C(v'_r)\in V$. Note that $(v''_r)$ is a 
basis of $V$ and
$$C\che{C}v''_r=C\che{C}C(v'_r)=C((-1)^{n+1}c_r\i v'_r)=(-1)^{n+1}c_r\i v''_r.$$
Thus $v''_r$ is an eigenvector of $C\che{C}$ that is $v''_r=\x_rv_{\t(r)}$ for some
$\x_r\in\RR-\{0\}$ and some permutation $\t$ of $[1,n]$ such that $c_{\t(r)}=c_r\i$. Since 
$c_1>c_2>\do>c_n>0$ it follows that $c_{\t(1)}<c_{\t(2)}<\do<c_{\t(n)}$ hence $\t(r)=r^*$
and $c_{r^*}=c_r\i$. Thus we have $C(v'_r)=\x_rv_{r^*}$ so that 
$\che{C}v_r=\x_r\i v'_{r^*}$ and $C\che{C}v_r=\x_{r^*}\x_r\i v_r$. We see that 
$\x_{r^*}\x_r\i=(-1)^{n+1}c_r$. We have 
$$\la v_r,v_s\ra=C\i(v_r)(v_s)=\x_{r^*}\i v'_{r^*}(v_s)=\x_{r^*}\i\d_{r,s^*}$$
and (a) follows.

\subhead 3\endsubhead
It is natural to ask whether the definitions of $G_{\ge0},G_{>0}$ given in \S1 are specific
to $GL(V)$ or have a Lie theoretical meaning. Thus we wish to replace $G$ in \S1 by an 
arbitrary split reductive connected algebraic group $G$ defined over $\RR$, such as a 
general linear group, a symplectic group or a group of type $E_8$, with a fixed pinning. 
(We shall identify an algebraic variety defined over $\RR$ with its set of $\RR$-rational 
points. Then $G$ is a not necessarily connected Lie group.) The pinning is a substitute for
the choice of a basis for $V$ in \S1; it is a family of homomorphisms $x_i,y_i$ from $\RR$
into $G$ indexed by a finite set $I$. (For a precise definition of a pinning see
\cite{\LU, \S1}; here we note only that for $G=GL(V)$, $x_i,y_i$ may be taken as in \S1.) 
Let $U^+$ (resp. $U^-$) be the subgroup of $G$ generated by $x_i(a)$ (resp. by $y_i(a)$)  
for various $i\in I,a\in\RR$. Let $T$ be the unique (split) maximal torus of $G$ which 
normalizes both $U^+$ and $U^-$. Let $B^+=TU^+$, $B^-=TU^-$ so that $T=B^+\cap B^-$. Let 
$T_{>0}$ be the "identity component" of $T$ (when we say "identity component" we mean this
in the sense of Lie groups, not algebraic groups). This agrees with the notation in \S1 for
$G=GL(V)$.

The results of Ann Whitney for $GL(V)$ suggest a definition in the general case for the 
monoid $G_{\ge0}$ (by repeating 1(b)) and for $G_{>0}$ as follows (see \cite{\LU}). Let 
$N=\dim U^+=\dim U^-$. We can find a sequence $\ii=(i_1,i_2,\do,i_N)$ in $I$ such that the
map $\RR_{>0}^N@>>>U^+$, 
$$(a_1,a_2,\do,a_N)\m x_{i_1}(a_1)x_{i_2}(a_2)\do x_{i_N}(a_N)$$ 
is injective; the image of this map is independent of the choice of $\ii$. It is an open 
submonoid of $U^+$ denoted by
$U^+_{>0}$. We define an open submonoid $U^-_{>0}$ of $U^-$ in a similar way (replacing 
$x_i$ by $y_i$). Let $G_{>0}=U^+_{>0}T_{>0}U^-_{>0}$. We have also 
$G_{>0}=U^-_{>0}T_{>0}U^+_{>0}$ and $G_{>0}$ is an open submonoid of $G$. We see that 
$G_{\ge0},G_{>0}$ do indeed have a meaning for an arbitrary $G$. As in \cite{\LU} we set 
(imitating the definition for $GL(V)$ given in \cite{\GK}) 
$G_{\ge0}^{\text{osc}}=\{g\in G_{\ge0};g^m\in G_{>0}\text{ for some }m\ge1\}$. We have
$G_{>0}\sub G_{\ge0}^{\text{osc}}\sub G_{\ge0}$ and $G_{\ge0}^{\text{osc}}$ is closed under
multiplication (see \cite{\LU, 2.19}).

The following analogue of the Gantmacher-Krein theorem for a general $G$ is contained in
\cite{\LU, 5.6}, \cite{\LU, 8.10}.

(a) {\it Let $g\in G_{>0}$. Then $g$ is contained in a unique $G$-conjugate of $T_{>0}$.}
\nl
We will now make some comments on the proof of (a). The key case in the proof is that where
$G$ is simply connected as an algebraic group and "simply-laced" (the last condition means
that $G$ is a product of groups of type $A,D,E$ in the Cartan-Killing classification); 
there are standard techniques by which various statements for a general $G$ can be reduced
to this case.

As in the proof of the classical Gantmacher-Krein theorem which was based on Perron's 
theorem we find that we need a definition of $G_{>0}$ along the lines of the original 
definition of Schoenberg in terms of minors of a matrix. In the case of $GL(V)$ that 
definition exploits the fact that the basic representations of $GL(V)$ (the exterior powers
of the standard representation) have a simple natural basis. For general (simply connected,
simply laced) $G$ we replace these exterior power representations by the finite dimensional
irreducible algebraic representations of $G$. Quite surprisingly it turns out that these 
representations admit canonical bases with respect to which the elements $x_i(a),y_i(a)$ 
with $a\ge0$ act by matrices with all entries in $\RR_{\ge0}$. This allows us to give a 
Schoenberg-style definition of $G_{\ge0}$ and $G_{>0}$ using these entries instead of 
minors. Then Perron's theorem is applicable for each of these representations and (a) 
follows. 

Note that there is no known definition of the canonical bases which is purely in terms of 
$G$. The only known definition \cite{\CB} uses the fact that the irreducible
representations of $G$ considered above are limits as $q\to1$ of irreducible 
representations of a "quantum group" which are some entities depending on a parameter $q$.
The canonical bases are first defined at the level of quantum group and then one takes
their limit as $q\to1$ to obtain the canonical bases for representations of $G$. Moreover
the positivity properties of the action of $x_i(a),y_i(a)$ with respect to the canonical 
basis come from a stronger property which holds for the generators of the quantum group. 
This stronger property is established using a geometric interpretation of the quantum group
and the associated canonical basis in terms of the theory of perverse sheaves and it 
ultimately depends on the theory around the Weil conjectures for algebraic varieties over a
finite field (proved by Deligne). Thus the statement (a) which is elementary in the case of
$GL(V)$ is very far from elementary in the general case.

\subhead 4\endsubhead
Let $\cb$ the set of all Borel subgroups of $G$ that is subgroups that are conjugate to 
$B^+$ (or equivalently to $B^-$). This is a compact manifold with a transitive $G$-action 
(conjugation) called the flag manifold. In \cite{\LU, \S8} the positive part $\cb_{>0}$ is
defined. It is the set of all $B\in\cb$ such that $B=uB^+u\i$ for some $u\in U^-_{>0}$ or 
equivalently such that $B=u'B^-u'{}\i$ for some $u\in U^+_{>0}$. In the case where 
$G=GL(V)$, $\dim V=2$, $\cb$ is a circle with two distinguished points $B^+,B^-$ and 
$\cb_{>0}$ is one of the two connected components of the complement of $\{B^+,B^-\}$ (a 
half circle). In general, $\cb_{>0}$ is an open ball in $\cb$. 

For $i\in I,a\in\RR$ we set ${}'x_i(a)=x_i(-a)$, ${}'y_i(a)=y_i(-a)$. Note that 
${}'x_i,{}'y_i$ define a new pinning for $G$. The objects attached to this new pinning in 
the same way as
$$B^+,B^-,U^+,U^-,T,T_{>0},U^+_{>0},U^-_{>0},G_{>0},\cb_{>0}$$
were attached to $x_i,y_i$ are:
$$B^+,B^-,U^+,U^-,T,T_{>0},{}'U^+_{>0}=(U^+_{>0})\i,{}'U^-_{>0}=(U^-_{>0})\i,
{}'G_{>0}=G_{>0}\i,{}'\cb_{>0}.$$
Recall that $B,B'$ in $\cb$ are said to be opposed if their intersection is a maximal 
torus. We show:

(a) {\it If $B\in\cb_{>0}$, $B'\in{}'\cb_{>0}$, then $B,B'$ are opposed.}
\nl
We have $B=uB^+u\i$ where $u\in U^-_{>0}$ and $B'=u'B^+u'{}\i$ where $u'\in{}'U^-_{>0}$. It
is enough to show that $u'{}\i Bu'=\tu B^+\tu\i$ is opposed to $B^+$ where $\tu=u'{}\i u$.
Since $u,u'{}\i$ are in $U^-_{>0}$ and $U^-_{>0}$ is closed under multiplication, we have 
$\tu\in U^-_{>0}$. Therefore $\tu B^+\tu\i,B^+$ are opposed by \cite{\LU, 2.13(a)}.

\subhead 5\endsubhead
In this section we describe a strengthening of 3(a) which also places 2(a) in a Lie 
theoretic setting.  

In the setup of \S3 let $\s:G@>>>G$ be an automorphism of finite order (as an algebraic 
group) which preserves the pinning that is, for some permutation $\t:I@>>>I$ we have 
$\s(x_i(a))=x_{\t(i)}(a)$, $\s(y_i(a))=y_{\t(i)}(a)$ for all $i\in I,a\in\RR$. We have the
following strengthening of 3(a).

(a) {\it Let $g\in G_{\ge0}^{\text{osc}}$. Define $\a:G@>>>G$ by $h\m g\s(h)g\i$. There is
a unique $B\in\cb_{>0}$ and a unique $B'\in{}'\cb_{>0}$ such that $B$ and $B'$ are 
$\a$-stable. Then $B,B'$ are opposed. Moreover, $\a$ induces a dilation on 
$\Lie(B)/\Lie(B\cap B')$, a contraction on $\Lie(B')/\Lie(B\cap B')$ and an automorphism 
of finite order of $\Lie(B\cap B')$.}
\nl
(An endomorphism $A:E@>>>E$ of a finite dimensional $\RR$-vector space is said to be a
dilation (resp. contraction) if all its eigenvalues $\l$ satisfy $|\l|>1$ (resp. 
$|\l|<1$).)

Let $m\ge1$ be such that $\s^m=1$. Let $g'=g\s(g)\s^2(g)\do\s^{m-1}(g)$. From the 
definitions we see that $\s(G_{\ge0})=G_{\ge0}$, $\s(G_{>0})=G_{>0}$,
$\s(G_{\ge0}^{\text{osc}})=G_{\ge0}^{\text{osc}}$. Since $G_{\ge0}^{\text{osc}}$ is closed
under multiplication it follows that $g'\in G_{\ge0}^{\text{osc}}$. Replacing $m$ by a 
multiple we can assume that we have $g'\in G_{>0}$.

Define uniquely $B\in\cb_{>0}$ by $g'\in B$, see \cite{\LU, 8.9}; define uniquely 
$B'\in{}'\cb_{>0}$ by $g'{}\i\in B'$, see \cite{\LU, 8.9} for ${}'G_{>0}$ instead of 
$G_{>0}$. Let $B_1=g\s(B)g\i\in\cb$, $B'_1=\s\i(g\i B'g)\in\cb$. From the definitions we 
see that $\s(\cb_{>0})=\cb_{>0}$, $\s({}'\cb_{>0})={}'\cb_{>0}$. Since 
$g\in G_{\ge0}$ we have $\Ad(g)(\cb_{>0})\sub\cb_{>0}$, see \cite{\LU, 8.12}.
Similarly since $g\i\in{}'G_{\ge0}$ we have $\Ad(g\i)({}'\cb_{>0})\sub{}'\cb_{>0}$.
We see that $B_1\in\cb_{>0}$, $B'_1\in{}'\cb_{>0}$. Applying $\s$ to $g'\in B$ and $\s\i$ 
to $g'{}\i\in B'$ we obtain $\s(g')\in\s(B)$, $\s\i(g'{}\i)\in\s\i(B')$. Since 
$\s(g')=g\i g'g$, $\s\i(g'{}\i)=\s\i(g)g'{}\i\s\i(g\i)$ we see that $g\i g'g\in\s(B)$, 
$\s\i(g)g'{}\i\s\i(g\i)\in\s\i(B')$ that is $g'\in B_1$, $g'{}\i\in B'_1$. By the 
uniqueness of $B,B'$ it follows that $B_1=B,B'_1=B'$. Thus we have $B=\a(B)$, 
$B'=\a\i(B')$. Hence $B'=\a(B')$. Note that $B\in\cb_{>0}$ is uniquely determined by the 
condition $B=\a(B)$. Indeed this condition implies that $B=\a^m(B)$ that is $g'\in B$ and 
this condition is known to determine $B$ uniquely. Similary $B'\in{}'\cb_{>0}$ is uniquely
determined by the condition $B'=\a(B')$. Thus the first assertion of (a) is established. 
The second assertion of (a) follows from 4(a). The third assertion follows from the 
analogous assertion where $\a$ is replaced by $\a^m=\Ad(g')$; hence to prove it we may 
assume that $m=1$ and $g\in G_{>0}$.
In this case, by \cite{\LU, 8.10} we can find unique $u\in U^-_{>0}$, $\tu\in U^+_{>0}$, 
$t\in T_{>0}$ such that $g=u'\tu tu\i$; moreover, $\Ad(t)$ is a dilation on $\Lie(U^+)$ and
a contraction on $\Lie(U^-)$. It follows that $\tu t$ is conjugate to $t$ under an element
in $U^+$ hence $\Ad(\tu t)$ is a dilation on $\Lie(U^+)$. We have $uB^+u\i\in\cb_{>0}$, 
$g\in uB^+u\i$ hence $uB^+u\i=B$ and $u\i gu=\tu t$. It follows that $\Ad(g)$ is a dilation
on $\Lie(uU^+u\i)=\Lie(B)/\Lie(B\cap B')$. Since $B,B'$ are opposed and $g\in B\cap B'$, 
the linear map $\Ad(g)$ on $\Lie(B')/\Lie(B\cap B')$ may be identified with the transpose 
inverse of $\Ad(g)$ on $\Lie(B)/\Lie(B\cap B')$ hence is a contraction. Since 
$g\in B\cap B'$, $\Ad(g)$ acts trivially on $\Lie(B\cap B')$. This completes the proof of
(a).

We now show:

(b) {\it Let $g\in G_{>0}$. Let $B,B'$ be as in (a). Then $g$ belongs to the
"identity component" of the torus $B\cap B'$.}
\nl
Let $u,\tu,t$ be as in the proof of (a). Since $u\i gu=\tu t$, $\tu t$ is contained in the
maximal torus $B^+\cap u\i B'u$ of $B^+$ and it is enough to show that it is contained in 
the "identity component" of $B^+\cap u\i B'u$ or equivalently its image in $B^+/U^+$ is 
contained in the "identity component" of $B^+/U^+$. But that image is the same as the image
of $t$ and it remains to use the fact that $t\in T_{>0}$.

\subhead 6\endsubhead
We now assume that $G$ is adjoint (that is with trivial centre) as an algebraic group. Let
$S$ be a closed Riemann surface of genus $g\ge2$. Let $\p_1$ be the fundamental group of 
$S$ at some point of $S$. In this section we review some recent results of Fock and 
Goncharov \cite{\FG} which show that total positivity can be used to understand some 
features of the real algebraic variety $\Hom(\p_1,G)$ of homomorphisms $\p_1@>>>G$.

Let $X$ be a $4$-elements set. A dihedral order on $X$ is a partition of $X$ into two
$2$-element sets. A map $f:X@>>>\cb$ is said to be {\it positive} (see \cite{\FG, \S5} if 
for some/any numbering $x_1,x_2,x_3,x_4$ of $X$ such that the dihedral order is 
$\{x_1,x_3\},\{x_2,x_4\}$ and some $h\in G$ we have
$$hf(x_1)h\i=B^+,hf(x_3)h\i=B^-,hf(x_2)h\i\in\cb_{>0},hf(x_4)h\i\in{}'\cb_{>0}.$$
(The equivalence of "some" and "any" is proved in \cite{\FG, Thm. 5.3}.)

Let $C$ be the unit circle in $\CC$. A map $F:C@>>>\cb$ is said to be a {\it positive 
curve} (see
\cite{\FG, Def. 6.4} if for any $4$-element subset $X$ of $C$ the restriction of $F$ to $X$
is positive with respect to the dihedral order on $X$ given by the partition $X=X'\sqc X''$
such that the line spanned by $X'$ meets the line spanned by $X''$ in a point of the (open)
unit disk. Let $\cc(\cb)$ be the set of all positive curves $C@>>>\cb$.

Now there is a free holomorphic action of $\p_1$ on the open unit disk whose orbit space is
$S$; this extends continuously to an action of $\p_1$ on the boundary $C$ of the open unit 
disk. This last action is such that for any $\g\in\p_1$ and any $4$-element subset $X$ of
$C$ the bijection $X@>>>\g(X)$ induced by $\g:C@>>>C$ is compatible with the dihedral order
on $X$ and $\g(X)$ defined as above. Hence if $F:C@>>>\cb$ is a positive curve then 
$F\circ\g:C@>>>\cb$ is a positive curve. Thus $\g:F\m F\circ\g\i$ is an action of $\p_1$ on
$\cc(\cb)$. On the other hand $G$ acts on $\cc(\cb)$ by $h:F\m{}^hF$ where
${}^hF(z)=hF(z)h\i$ for any $z\in C$. This $G$-action commutes with the $\p_1$-action and
is free: if $h\in G$, $F\in\cc(\cb)$ satisfy ${}^hF=F$ then $h$ is contained in three Borel
subgroup any two of which are opposed; hence $h=1$. Let $\cc_{\p_1}(\cb)$ be the set of all
$F\in\cc(\cb)$ such that for any $\g\in\p_1$, $F\circ\g$ is in the $G$-orbit of $F$ that 
is, $F\circ\g={}^{\c_F(\g)}F$ for some $\c_F(\g)\in G$ which is unique (by the freeness of 
the $G$-action). Note that $\c_F:\p_1@>>>G$ is a homomorphism. The map 
$\cc_{\p_1}(\cb)@>>>\Hom(\p_1,G)$, $F\m\c_F$ is compatible with the $G$-actions (where $G$
acts on $\Hom(\p_1,G)$ by conjugation). Let $\Hom^{\text{pos}}(\p_1,G)$ be the image of
this map (the set of "positive homomorphisms"). By passage to $G$-orbits we get a map
$\Ph:G\bsl\cc_{\p_1}(\cb)@>>>G\bsl\Hom^{\text{pos}}(\p_1,G)$.

According to \cite{\FG}, $\Ph$ is a bijection and $G\bsl\Hom^{\text{pos}}(\p_1,G)$ is a
ball of dimension $(2g-2)\dim G$; moreover, any positive homomorphism $\c:\p_1@>>>G$ is 
injective with discrete image and for any $\g\in\p_i-\{1\}$, $\c(\g)$ is contained in a
$G$-conjugate of $G_{>0}$ hence, by 3(a), it is contained in a unique $G$-conjugate of 
$T_{>0}$.

This sheds a new light on a result of Hitchin \cite{\HI}. Let $\Hom^{cr}(\p_1,G)$ be the 
space of all $\c\in\Hom(\p_1,G)$ such that the induced action of $\p_1$ on $\Lie(G)$ is 
completely reducible. In \cite{\HI}, Hitchin shows using techniques of analysis that there 
is a canonical connected component of $G\bsl\Hom^{cr}(\p_1,G)$ which is a ball of dimension
$(2g-2)\dim G$. As a consequence of \cite{\FG}, this "Hitchin component"
coincides with $G\bsl\Hom^{\text{pos}}(\p_1,G)$.

\subhead 7\endsubhead
Let $P^m$ be the $m$-dimensional real projective space. Following Schoenberg we say that an
imbedding $f:C@>>>P^m$ is a {\it convex curve} if for any hyperplane $H$ in $P^m$ the 
intersection $f(C)\cap H$ has at most $n$ points. In \cite{\FG, Thm. 9.4} it is shown that
if $G=GL(\RR^{m+1})$ then the image of a positive curve in $\cb$ under the natural map
$\cb@>>>P^m$ is a convex curve; moreover, with suitable smoothness assumptions, this gives
a bijection between positive curves in $\cb$ and convex curves in $P^m$. 

{\it Acknowledgement.} I thank A. Goncharov for useful discussions.

\enddocument